\documentclass[oneside,reqno]{amsart}
\usepackage[english]{babel}
\usepackage{indentfirst}
\usepackage[T1]{fontenc}
\usepackage[utf8]{inputenc}
\usepackage{paralist}
\usepackage{amssymb}
\usepackage{amsthm}
\usepackage{amsmath}
\usepackage{amscd}
\usepackage{graphicx}
\usepackage[colorlinks=true]{hyperref}
\usepackage[margin=1.5in]{geometry}

\theoremstyle{plain} 
\newtheorem{teo}{Theorem}[section]

\newtheorem{prop}[teo]{Proposition}
\newtheorem{defn}[teo]{Definition}

\numberwithin{equation}{section} 

 \renewcommand{\phi}{\ensuremath{\varphi}}

\begin{document}
\nocite*
\date{}

\author{Fausto Ferrari}

\address{Fausto Ferrari: Dipartimento di Matematica\\ Universit\`a di Bologna\\ Piazza di Porta S.Donato 5\\ 40126, Bologna-Italy}
\email{fausto.ferrari@unibo.it }
\thanks{The author wishes to thank Salvatore Coen and Salomon Ofman for sharing with him their information about the life of French mathematicians between the two world wars of the 20-th century and Bruno Franchi for telling him the role of some public officers in French University system. Moreover the author also thanks Catherine Goldstein for pointing out the thesis  \cite{Leloup} to him and for having told him her opinion about some questions concerning the life of A. P. Marchaud}

\thanks{The author is supported by  MURST, Italy, and
GNAMPA project 2017: {\it Regolarit\`a delle soluzioni viscose per equazioni a derivate parziali non lineari degeneri}}

\title[]{Weyl and Marchaud derivatives: a forgotten history}
\begin{abstract}
In this paper we recall the contribution given by Hermann Weyl and Andr\'e Marchaud to the notion of fractional derivative. In addition we discuss some relationships between the fractional Laplace operator and Marchaud derivative in the perspective to generalize these objects to different fields of the mathematics.
\end{abstract}

\maketitle

\tableofcontents

\section{Introduction}
Exactly one century ago, while we are writing, in 1917, a paper by Hermann Weyl, \emph{Bemerkungen zum Begriff des Differentialquotienten gebrochener Ordnung}, \cite{Weyl},  appeared. It dealt with the definition of a fractional derivative in a weaker sense with respect to the approach classically known at that time with the name of Riemann-Liouville derivative.

Ten years later, in 1927, it was published the thesis of a misunderstood French mathematician, Adr\'e Paul Weyl who discussed at the age of forty his PhD work entitled \emph{ Sur les d\'{e}riv\'ees et sur les diff\'{e}rences des fonctions de
  variables r\'{e}elles}, \cite{Marchaud_T}.
  
In \cite{SKM} the names Weyl and Marchaud  appear associated with the notion of fractional derivative more than two hundred times. Nevertheless, in my opinion,  the name Marchaud is not so popular even among the mathematicians dealing with fractional calculus, in particular among scientists coming from western countries. Due to the huge quantity of papers dealing with fractional subjects, my previous statement could appear debatable. In any case this opinion can be tested just consulting, for instance a database. We tried, for instance, with the AMS database MathShiNet.  In fact, inserting the keyword "Marchaud" anywhere, we obtain around two hundred files.  Among these two hundred files, improving the request by searching also the word Marchaud in the titles of the papers, we find around fifty files. In addition, by reading these titles, covering for example the last twenty years, we realize at a first glance that the frequency of mathematicians from Eastern countries is prevalent. Indeed, on the contrary of what we stated about western mathematicians, Marchaud's name is recurrent in fractional calculus literature and among mathematicians coming from Eastern Europe, let us recall one more time the number of citations that appear in \cite{SKM}.  

Concerning Hermann Weyl,  of course, we are considering a very popular mathematician for many other mathematical reasons. Nevertheless, we have to say that also in this case Weyl's name is not usually associated with the fractional calculus even if the specialists in the field are aware of the importance of his contribution in fractional calculus. It could be interesting to understand whether many of them know why a fractional derivative is entitled to him, but this is another story.
  
 For different reasons the authors of the two cited papers will not publish any more results explicitly amenable to the fractional derivative. So, accepting the previous interpretation, they appear as isolated points in the {\it mare magnum} of the fractional calculus, where the more popular names are nowadays others.
   
In Marchaud doctoral thesis, see p. 47, Section 27 formula (23) in \cite{Marchaud_T}, or the definition (23) in the published paper \cite{Marchaud} at p. 383, he defined the following fractional differentiation for sufficiently regular real functions $f:(0,1)\to \mathbb{R}$ extended with $0$ for $x\leq 0,$ whenever $\alpha\in (0,1):$
$$
{\bold D}^{\alpha}f(x)=\frac{\alpha}{\Gamma(1-\alpha)}\int_0^{+\infty}\frac{f(x)-f(x-t)}{t^{1+\alpha}}dt.
$$
This definition can be easily given for a function defined in all of $\mathbb{R}$ and for every $\alpha\in (0,1)$ distinguishing two types of derivative, see \cite{SKM}, respectively from the right  and  from the left: 
$$
{\bold D}^{\alpha}_+f(x)=\frac{\alpha}{\Gamma(1-\alpha)}\int_0^{+\infty}\frac{f(x)-f(x-t)}{t^{1+\alpha}}dt
$$
and
$$
{\bold D}^{\alpha}_-f(x)=\frac{\alpha}{\Gamma(1-\alpha)}\int_0^{+\infty}\frac{f(x)-f(x+t)}{t^{1+\alpha}}dt.
$$
The construction of these operators will be briefly described in the next section following the original motivation contained in Marchaud's thesis. 

The problem of giving a coherent definition of derivative of a function for all positive real numbers has a long history. For instance see \cite{Ross}, \cite{SKM}, \cite{KST06}, \cite{Butzer_Westphal}, \cite{BDST} for some detailed information. In any case Abel's contribution for solving the tautochrone problem, \cite{Abel}, and the work by Liouville \cite{Liouville_J} and Riemann \cite{Riemann_B} in application to geometry are fundamental and well known at the beginning of the fractional calculus. Many other authors have written papers that contributed to improve the knowledge of this subject. Nevertheless I think that a very special role has to be recognized to Hermann Weyl because, probably following the path traced by Riemann, as Weyl himself writes in \cite{Weyl}, he introduced, maybe first, the nonlocal operator that is known as Marchaud derivative, for people who know it, in a significative, even if particular, case. We shall dedicate Section \ref{Weyl_approach}  to this aspect.
 
 Our interest to this subject comes out after the celebrated contribution given by the paper \cite{CS}. Indeed, the authors developed an idea that was already contained in \cite{Molc_Ostrov}. In any case, in \cite{CS}, the authors dedicated their interest to a different type of nonlocal operator with respect to the fractional derivative: the fractional Laplace operator. In particular in \cite{CS} a different perspective in the interpretation of the nonlocal operators was introduced using a method based on an {\it extension} approach, see also \cite{Molc_Ostrov}. We do not want to bore the reader too much with this subject. However, some words are in order. Heuristically, following the extension approach idea it is possible to deduce the properties of a nonlocal operator from the ones of a local operator. In \cite{CS} the authors were concerned with the fractional Laplace operator, while the local operator obtained after the extension construction was a degenerate elliptic operator in divergence form. This approach can be developed considering the solution of {\it ad hoc} Dirichlet problem formulated in an unbounded set, where an auxiliary variable has been added, and then taking the limit of a weighted normal derivative of  the solution of the Dirichlet problem, when this auxiliary variable vanishes. The scientific follow up of \cite{CS} produced an enormous amount of papers. Moreover, in \cite{StiTor} such idea was generalized considering an abstract approach in a very powerful way. Following this stream of ideas,  in \cite{FerFra},  an intrinsic characterization of the fractional sub-Laplace operators in Carnot groups was obtained. Roughly saying, the operators considered in this last case are sums of squares of smooth vector fields satisfying the H\"ormander condition in a non-commutative structure. 
 
 The approach described in \cite{CS} and \cite{StiTor}, and then in \cite{FerFra}, was also extended to the case of fractional operators in \cite{Torreault} and independently also, as very often it happens when the time is ripe, in \cite{BF}.  
 Indeed, with this aim, commenting for instance \cite{BF}, we faced the problem of defining the Marchaud derivative via an extension approach in order to obtain a Harnack inequality for solutions of homogeneous fractional equations.  As a consequence of this research we realized in particular that Marchaud derivative and Weyl derivative have been, in a sense, perhaps a little put aside in the last time, especially considering the great development and the large popularity that researches about nonlocal operators have recently had.
 This last remark is essentially based on the popularity of other fractional derivatives, for instance the Riemann-Liouville  derivative or the Caputo derivative, see \cite{LG}, \cite{KST06} and \cite{GLY}, for a modern approach to these operators. Indeed, see also \cite{ACV} for a recent example of application involving Caputo derivative.
  
 On the other hand, by reading the monumental opera \cite{SKM},  it is possible to verify, as we pointed out at the beginning of this introduction,  that  Weyl and Marchaud names are cited many times. So, the curiosity of explaining this situation was strong. Why  do only few people associate Weyl  and Marchaud names to the fractional subject? More precisely, why do only few people utilize these fractional derivatives for applications, simply preferring other definitions, even if it appears natural to use Weyl and Marchaud operators?  We do not have any conclusive reply. In any case, it is quite difficult to understand what are the true motivations of this apparent amnesia. Of course the specialists of the fractional calculus know Marchaud and Weyl derivatives in reading in particular \cite{SKM} where the right tribute to both these mathematicians has been given. Perhaps only recently a new awaken interested about these definitions spread out. For this reason we think,  hopefully, that the contribution of this paper might be useful in consolidating this new trend. 
 
Anyhow, a partial reply to the previous questions, partially related with beginning of these events, can be found in the social contour that strongly influenced the life of the two mathematicians.
 The period during which these researches were developed was very uproarious for Europe. Weyl's paper \cite{Weyl} dated back to 1917, Marchaud's thesis \cite{Marchaud_T} was published in 1927 and both the lives of these two people were, from different reasons, effected by the two world wars events, see e.g. \cite{EnciWeyl} for some biography details about Weyl's life and \cite{Condette} and \cite{Leloup} for some information about Marchaud. 

In this paper we want to analyze the definition of fractional derivative given by Weyl and Marchaud, concentrating on those  aspects  that, in perspective, seem to be more flexible for generalizing to other situations the notion of nonlocal operator, see e.g. \cite{BalaAV} for facing the case of the semigroup approach in its abstract generality and then for recalling the contribution given in \cite{StiTor}, for fractional Laplace operators, and then recalling \cite{GMS} for several generalizations in an abstract approach, including in principle: Riemannian manifolds, Lie groups, infinite dimensions and non symmetric operators, see also \cite{FerFra} for the particular case of Carnot groups. 
We would also like to point out  some recent researches establishing few relationships between Marchaud derivative and some nonlocal operators in non-commutative structures, see  \cite{Fer_H}. 

With this aim we summarize the plan of the paper. After this introduction, the reader can find in Section \ref{bioWM} some basic biographic information about Weyl and Marchaud. In Section \ref{Marchaud_approach} we discuss briefly how Marchaud came to define his derivative following by reviewing a part of his PhD thesis. In Section \ref{Grunwald_Letnikove_approach} we recall the basic idea already developed earlier by Gr\"unwald and Letnikov that is at the base of the fractional derivative given by Weyl and Marchaud. In Section \ref{Weyl_approach} we recall the seminal Weyl's paper \cite{Weyl} discussing some details about the relationship between  his contribution and Marchaud's derivative. In Section \ref{ideasMW} we comment the basic ideas of the respective definitions. We face the modern general setting of  Marchaud derivative in Section \ref{general_setting} making also some remarks about its properties with respect to PDEs. In Section  \ref{fractional_recall} we continue our work by recalling the definition of fractional Laplace operator, while in Section \ref{FerBucSec} we deal with the definition of Marchaud derivative via an extension approach and eventually, in Section \ref{rel_Mar_Lapl}, we conclude our effort putting in evidence the relationship between the fractional Laplace operator and the Marchaud derivative.

In order to outline few aspects of Weyl's and Marchaud's biographies we list below only some key facts of the period during which the results about fractional derivative were written. For further curiosities or remarks we suggest to consult \cite{EnciWeyl}, \cite{Condette} and \cite{Leloup}. 

Closing this introduction we remark that from a historical point of view it would be interesting to deepen our knowledge of these two characters of the mathematical world, especially considering the influence and the role of the respective mathematical schools compared to the other mathematicians of their time and their scientific legacy.

\section{Short historical placement}\label{bioWM}

In this section we introduce some information about the lives of Weyl and Marchaud, mainly regarding the period of publication of their papers on fractional derivative without pretending to consider this parallel description exhaustive.

Hermann Weyl was born in Germany in 1885. Andr\'e Paul Marchaud was born in France in 1887. 

In 1913 Weyl was professor at the ETH (Swiss Federal Institute of Technology) in Z\"urich where he interacted also with Einstein. In 1915 Weyl was called up for military service in Germany, but since 1916 he was exempted from military duties for reasons of health. Later on he came back to Germany as successor of Hilbert in G\"ottingen, but in 1933 he left to go to IAS in Princeton, escaping from the Nazi regime, where he continued his brilliant career, see e.g. \cite{EnciWeyl} until he died in Z\"urich in 1955.

In 1913, Marchaud had not got his PhD thesis yet, probably because of his health problems. He was professor in a lyceum when he was mobilized by the French army in 1914. In the same year Marchaud was taken prisoner. He stayed in a {\it Oflag}  (a prison camp for officers only) from 1914 to 1918, at the beginning in Germany and then,  thanks to the help of the Red Cross who intervened because he was ill,  since 1917 in Switzerland. Marchaud discussed his PhD thesis later on, only in 1927. He continued his career mainly serving as {\it Rector} (provost) of French universities, even during the Nazi occupation of France in the Second World War,  until 1957, see \cite{Condette}, when he retired. He died in Paris in 1973.

\section{The Marchaud approach}\label{Marchaud_approach}

As we have already announced, in this section we represent the Marchaud approach following the main steps of a part of his PhD thesis.

The Liouville-Riemann integral of order $\alpha>0$ of a function $f:[a,b]\to \mathbb{R}$ is defined as:
$$
I^{(\alpha)}_af(x)=\frac{1}{\Gamma(\alpha)}\int_0^{x-a}t^{\alpha-1}f(x-t)dt.
$$
In this case the derivative of order $\alpha<n$ where $n$ is a positive integer is defined by
$$
\mathcal{D}^{\alpha}_a f(x)=D^nI_a^{n-\alpha}f(x).
$$
In fact this definition is well posed because it is independent to $n.$ In particular it is coherent since 
the following fundamental identity holds for all the functions in $L^1([a,b])$  that are bounded, for every $\alpha, \alpha'>0:$
$$
I_a^{(\alpha)}[I_a^{(\alpha')}f(x)]=I^{(\alpha+\alpha')}_af(x).
$$ 
The point of starting by Liouville, as Marchaud observed in \cite{Marchaud}, is that if $a=-\infty,$ then
$$
\mathcal{D}^{\alpha}_{-\infty}e^{kx}=k^\alpha e^{kx}.
$$
Moreover for $\beta\geq 0$ it results that
$$
I^{\beta}_a(1)=\frac{(x-a)^\beta}{\Gamma(\beta+1)},
$$
and in particular $I^{0}_a(1)=1,$
so that for every $\beta\geq 0,$ and for every $\alpha> 0$ we get
$$
\mathcal{D}^{\alpha}_aI^{\beta}_a(1)=D^{(n)}I_a^{n-\alpha}I^{\beta}_a(1)=D^{(n)}I_a^{n-\alpha+\beta}1=D^{(n)}\frac{(x-a)^{n-\alpha+\beta}}{\Gamma(n-\alpha+\beta+1)}=\frac{(x-a)^{\beta-\alpha}}{\Gamma(\beta-\alpha+1)}.
$$
Hence for $\beta=0$
$$
\mathcal{D}^{\alpha}_a(1)=\mathcal{D}^{\alpha}_aI^{0}_a(1)=\frac{(x-a)^{-\alpha}}{\Gamma(-\alpha+1)}.
$$
As a consequence the fractional derivative of order $\alpha$ for $\alpha$ non-integer is in general infinite for $x=a.$

Trying to define the fractional derivative as the fractional integral of negative order $\alpha$ we obtain a divergent integral.
In fact, formally, we should obtain
$$
I^{(-\alpha)}_af(x)=\frac{1}{\Gamma(-\alpha)}\int_0^{x-a}t^{-\alpha-1}f(x-t)dt.
$$
Then Marchaud argues in this way. Taking the integral of order $\alpha$ and assuming to consider the function extended with $0$ from $-\infty$ to $a,$ we get that 
$$
I^{(\alpha)}_{-\infty}f(x) \Gamma(-\alpha)=\int_0^{\infty}t^{\alpha-1}f(x-t)dt
$$
that is, making clear the definition of $\Gamma,$
$$
I^{(\alpha)}_{-\infty}f(x)\int_0^\infty t^{\alpha-1}e^{-t}dt=\int_0^{\infty}t^{\alpha-1}f(x-t)dt.
$$
The same formula holds for every positive integer $k,$ so that performing a change of variable like $t=ks$ in both the integrals we get: 
$$
I^{(\alpha)}_{-\infty}f(x) k\int_0^\infty (ks)^{\alpha-1}e^{-ks}ds=\int_0^{\infty}(ks)^{\alpha-1}f(x-ks)kds.
$$
that implies
$$
I^{(\alpha)}_{-\infty}f(x) \int_0^\infty s^{\alpha-1}e^{-ks}ds=\int_0^{\infty}s^{\alpha-1}f(x-ks)ds.
$$
Then taking a linear combination od order $p+1$ for a finite sequence of integer positive  decreasing number $\{k_i\}_{0\leq i\leq p}$we obtain summing terms by terms
$$
I^{(\alpha)}_{-\infty}f(x) \int_0^\infty s^{\alpha-1}\psi(s)ds=\int_0^{\infty}s^{\alpha-1}\phi(x,s)ds,
$$
where
$$
\psi(s)=\sum_{i=0}^{p}C_ie^{-k_is},\quad \phi(s)=\sum_{i=0}^{p}C_if(x-k_is),
$$
and $\{C_i\}_{1\leq 0\leq p}\subset \mathbb{R}.$
At this point Marchaud asks that passing to negative exponent $-\alpha$ the following relation makes sense
$$
I^{(-\alpha)}_{-\infty}f(x) \int_0^\infty s^{-\alpha-1}\psi(s)ds=\int_0^{\infty}s^{-\alpha-1}\phi(x,s)ds,
$$
calling the function $I^{(-\alpha)}_{-\infty}f(x)$ the fractional derivative of order $\alpha,$ that is $\mathcal{D}^{\alpha}f(x)$ is implicitly defined by
$$
\mathcal{D}^{\alpha}f(x) \int_0^\infty s^{-\alpha-1}\psi(s)ds=\int_0^{\infty}s^{-\alpha-1}\phi(x,s)ds.
$$
supposing that it is possible to chose $\psi$ in such a way $ \gamma(\alpha):=\int_0^\infty s^{-\alpha-1}\psi(s)ds$ does not vanish and as a consequence obtain the espression of $\phi.$ Discussing this problem
Marchaud find that if it is possible to find $\psi$ and $\phi$ with previous properties, then

$$
\gamma(\alpha)\mathcal{D}^{\alpha}f(x)=\int_0^{\infty}s^{-\alpha-1}\phi_\alpha(x,s)ds,
$$
where
$$
\phi_\alpha(x,s)=\sum_{i=0}^pC_{k_i}f(x-k_is).
$$
with a possible choice for $\psi$ that it is given by
$$
\psi(t)=e^{-t}(1-e^{-t})^p=\sum_{j=0}^p(-1)^{j}\binom{p}{j}e^{(-1-j)t}
$$
and
$$
\phi_\alpha(x,s)=\sum_{i=1}^p(-1)^{j-1}\binom{p}{j}f(x-js).
$$

After a  detailed computation Marchaud conclude that the existence of the fractional derivative $\mathcal{D}^{\alpha}f(x)$ continuous for continuous functions defined in $(a,b)$ is equivalent to the uniform convergence of the following integral 
$$
\int_\epsilon^{\infty}s^{-1-\alpha}\phi(x,s)ds
$$
 in every interval $(a',b)\subset (a,b)$ as $\epsilon\to 0^+$ and it is independent to the choice of the positive numbers $\{k_i\}_{1\leq i\leq p}.$
 A this point Marchaud define the fractional derivative of order $\alpha<p$ of a function defined in all of $\mathbb{R}$ implicitly:
 $$
 {\bold D}^{\alpha}f(x)\int_{0}^{\infty}s^{-1-\alpha}e^{-s}(1-e^{-s})^pds=\int_0^\infty s^{-\alpha-1}\sum_{j=1}^p(-1)^{j-1}\binom{p}{j}f(x-js)ds,
 $$
 or taking $\psi(t)=(1-e^{-t})^p-(1-e^{-2t})$ it is possible to obtain
 $$
 {\bold D}^{\alpha}f(x)\int_{0}^{\infty}s^{-1-\alpha}\left((1-e^{-s})^p-(1-e^{-2s})^p\right)ds=\int_0^\infty s^{-\alpha-1}(\Delta_{-s}^pf(x)-\Delta_{-2s}^pf(x))ds,
 $$
 where
 $$
\Delta_{-s}^pf(x)= \sum_{j=0}^p(-1)^{j}\binom{p}{j}f(x-(p-j)s).
 $$
 Hence, separating the integral,remarking that
  $$
\int_{0}^{\infty}s^{-1-\alpha}(1-e^{-2s})^pds=2^{\alpha}\int_{0}^{\infty}t^{-1-\alpha}(1-e^{-t})^pdt,
 $$
 and
 $$
 \int_0^\infty s^{-\alpha-1}\Delta_{-2s}^pf(x)ds= 2^{\alpha}\int_0^\infty t^{-\alpha-1}\Delta_{-t}^pf(x)dt
 $$
 we obtain
\begin{equation*}
\begin{split}
(1-2^\alpha) {\bold D}^{\alpha}f(x)\int_{0}^{\infty}s^{-1-\alpha}(1-e^{-s})^pds=(1-2^\alpha)\int_0^\infty s^{-\alpha-1}\Delta_{-s}^pf(x)ds,
 \end{split}
\end{equation*}
As a consequence we obtain also this representation
\begin{equation}\label{exampleMarchaud}
\begin{split}
{\bold D}^{\alpha}f(x)\int_{0}^{\infty}s^{-1-\alpha}(1-e^{-s})^pds=\int_0^\infty s^{-\alpha-1}\Delta_{-s}^pf(x)ds,
 \end{split}
\end{equation}

\section{Gr\"unwald-Letnikov derivative}\label{Grunwald_Letnikove_approach}
 It is impossible to deal with fractional Marchaud derivative without recalling the contribution of Gr\"unwald and Letnikov, see \cite{Grunwald}, and \cite{Letnikov}.   In fact for  giving a different perspective of the Marchaud derivative we have to introduce the Gr\"unwald-Letnikov derivative.\footnote{ Indeed, from this point of view, after we had completed this manuscript Francesco Mainardi pointed out the survey paper \cite{Rogosin_et_alt} dedicated to Marchaud and Gr\"unwald-Letnikov derivatives.} 
 To do this we need some new notation. 
 
 We recall that the binomial coefficients can be defined for every $\alpha\in \mathbb{C}$ and $n\in \mathbb{N}\cup\{0\}$ as:
\begin{equation}\label{altradef_iper}
\binom{\alpha}{0}=1,\quad \binom{\alpha}{n}=\frac{\alpha(\alpha-1)\cdots(\alpha-n+1)}{n!}=\frac{(-1)^n(-\alpha)_n}{n!},\quad n\in \mathbb{N}.
\end{equation}
It is also true that
$$
\binom{\alpha}{n}=\frac{\Gamma(\alpha+1)}{n!\Gamma(\alpha+n-1)}
$$
for $\alpha\in \mathbb{C}\setminus-\mathbb{N}$ and  $n\in \mathbb{N}.$

We introduce now the following notation concerning the difference of fractional order $\alpha\in \mathbb{R}$ for a function $f$ as follows. Let us denote
$$
(\Delta_h^\alpha f)(x)=\sum_{k=0}^{\infty}(-1)^k\binom{\alpha}{k}f(x-kh).
$$
We are now in position to define the Gr\"unwald-Letnikov fractional derivative, see \cite{Grunwald}, \cite{Letnikov} and also \cite{SKM}. Let $\alpha\in (0,1)$ be fixed and let $f:\mathbb{R}\to \mathbb{R}$ be a given function. The Gr\"unwald-Letnikov derivative of order $\alpha$ of $f$ is defined, separating the two cases, respectively as:
$$
f^{\alpha)}_+(x)=\lim_{h\to 0^+}\frac{(\Delta_h^\alpha f)(x)}{h^\alpha}
$$
and
$$
f^{\alpha)}_{-}(x)=\lim_{h\to 0^+}\frac{(\Delta^\alpha_{-h})f(x)}{h^\alpha},
$$
whenever the limit exists.

In order to understand better the reason of this definition we introduce the following definition.
\begin{defn}
We define a non-centered difference of increment $h$ on $f:\mathbb{R}\to \mathbb{R},$  as
$$
(I-\tau^{-t})f(x)=f(x)-f(x-t).
$$ 
\end{defn}
Then we obtain for every $m\in \mathbb{N}$ so that
$$
(I-\tau^{-t})^m=\sum_{k=0}^m(-1)^k\binom{m}{k}(\tau^{-t})^{k}
$$
and 
$$
(I-\tau^{-t})^mf(x)=\sum_{k=0}^m(-1)^k\binom{m}{k}(\tau^{-t})^{k}f(x)=\sum_{k=0}^m(-1)^k\binom{m}{k}f(x-kt).
$$
On the other hand, taking the Taylor expansion of the function $t\to (1+t)^\alpha$ in the center $t_0=0$ and $\alpha\in (0,1),$ we get
$$
(1+t)^\alpha=\sum_{k=0}^{+\infty}\binom{\alpha}{k}t^k,
$$
where
\begin{equation*}
\binom{\alpha}{0}=1,\quad \binom{\alpha}{n}=\frac{\alpha(\alpha-1)\cdots(\alpha-n+1)}{n!}=\frac{(-1)^n(-\alpha)_n}{n!},\quad n\in \mathbb{N}.
\end{equation*}
Thus, we can extend our definition to the fractional case it is possible to define
 for $\alpha\in(0,1)$
 $$
(I-\tau^{-t})^\alpha f(x)=\sum_{k=0}^{+\infty}(-1)^k\binom{\alpha}{k}(\tau^{-t})^{k}f(x)=\sum_{k=0}^{\infty}(-1)^k\binom{\alpha}{k}f(x-kt).
$$
In this way, we still maintain the semigroup property for the $\Delta_{h}^\alpha=(I-\tau^{-h})^\alpha,$ because for every $\alpha_1,\alpha_2\in \mathbb{R}$
$$
\Delta_{h}^{\alpha_1}\Delta_{h}^{\alpha_2}=(I-\tau^{-h})^{\alpha_1}(I-\tau^{-h})^{\alpha_2}=(I-\tau^{-h})^{\alpha_1+\alpha_2}
$$
and $(I-\tau^{-h})^{0}=I.$ Here we simply discuss the case of $\Delta_{+, h}^\alpha$ for $\alpha\in (0,1)$ but the results may be generalized to different exponents.

Moreover the following result holds, see e.g. \cite{KST06}.

\begin{teo}
Let $\alpha,\beta>0.$ Then for every bounded function:
$$
\Delta_h^\alpha\Delta_h^\beta f=\Delta_h^{\alpha+\beta}f.
$$
\end{teo}
In addition, considering one more time \cite{KST06}, and recalling also the contribution given in \cite{Chapman}, we have that:

\begin{teo}
Let $\alpha>0.$ Then, for every $f\in L^1(\mathbb{R})$
$$
\mathcal{F}(\Delta_h^\alpha f)(x)=(1-e^{ixh})^\alpha \mathcal{F}(f)(x).
$$
\end{teo}

In particular it is true that Gr\"unwald-Letnikov derivative of order $\alpha\in (0,1)$ coincides with Marchaud derivative of the same order. Indeed, in consideration of the two previous trivial properties the following result holds. 

\begin{teo}
Let $f\in L^p(\mathbb{R}),$ $p\geq 1.$ Then for every $q\geq 1$ there exist
$$
f^{\alpha)}_{\pm}(x)=\lim_{h\to 0,\:\mbox{in}\:\:L^q}\frac{\Delta^\alpha_{\pm h}f(x)}{h^\alpha}
$$
and
$$
\bold{D}^\alpha_{\pm}f(x)=\lim_{\epsilon\to 0,\:\mbox{in}\:\:L^q}C(\alpha)\int_{\epsilon}^{+\infty}\frac{f(x)-f(x\mp h)}{h^{1+\alpha}}dh
$$

Moreover 
$$
f^{\alpha)}_{\pm}(x)=\bold{D}^\alpha_{\pm}f(x),
$$
independently  to $p$ and $q.$
\end{teo}
The proof is quite long and can be find in \cite{SKM}, Theorem 20.4.  Moreover, about this topic we recall the very recent contribution \cite{ADT}. By the way, this last paper can be considered also as further signal of the renascent interest for Marchaud derivative. In fact in that manuscript has been recently proved  the coincidence of Marchaud derivative and Grunwald-Letnikov derivative for functions in H\"older spaces with explicit rates of convergence. 
Previous results encode many facts. The first concerns with the commutativity of Gr\"unwald-Letnikov derivative as well as the Marchaud  derivative, namely $(f^{\alpha)})^{\beta)}=(f^{\beta)})^{\alpha)}=f^{\alpha+\beta)}$ and ${\bold D}^\alpha{\bold D}^\beta={\bold D}^\beta{\bold D}^\alpha={\bf D}^{\alpha+\beta}.$ 
\section{Weyl Derivative}\label{Weyl_approach}

Hermann Weyl's name is associated with many important scientific results in physics and mathematics. In particular, concerning fractional derivative, Weyl gave an important contribution that is strictly linked to the Marchaud derivative. By the truth Weyl introduced in its paper \cite{Weyl}, at p. 302, exactly the definition of the fractional derivative that Marchaud gave in \cite{Marchaud}. The paper written by Weyl appeared in 1917, while Marchaud thesis was published in 1927, see \cite{Marchaud_T} and \cite{Marchaud}. It is not clear if the two definitions were discovered independently. The cited Weyl's paper, whose title is {\it Bemerkungen zum Begriff des Differentialquotienten gebrochener Ordnung} (Remarks on the notion of the differential quotient of a broken order), concerns with the notion of fractional derivative. In the introduction of his paper,  Weyl recognizes at first the efforts made by Bernanrd Riemann for obtaining  a notion of derivative for every positive real number. In particular Weyl cited the contain of the unpublished Riemann notes reported in the XIX paper of  the published Riemann opera post. About this fact, Weyl recalls as the editor of that volume remarks that Riemann surely did not think that those computations would have been published, at least in that form. In any case Weyl faces the problem of starting from those notes and having in mind that he wants to obtain a definition that works for periodic functions. In order to avoid the problem of introducing some privileged points, as very often happens in literature concerning fractional derivative, he assumes that periodic functions have to have zero mean. We do not enter into the details here (see Section ), however Weyl uses the properties of Fourier series and, at page 302, \cite{Weyl}, appears the following relationship
\begin{equation}\label{Weyl1}
g(x)=\beta\int_0^{\infty}\frac{f(x)-f(x-\xi)}{\xi^{1+\beta}}d\xi
\end{equation}
for functions $f$ H\"older having modulus of continuity with exponent $\alpha$ and $\alpha>\beta$ and knowing that $g$ denotes the fractional derivative of order $\beta.$ The same argument was reported in \cite{Zygmund} at page 226, see formula (3) in IX, 9.81. Nevertheless previous formula, apparently, disappeared in the final version of the book \cite{Zygmund_b} published later on, probably because the author was mainly interested in the periodic properties  of the functions, but we do not have any proof of this statement.

Anyhow, also in \cite{SKM} it can be find the definition of Weyl derivative. Starting from Fourier expansion of a periodic function, Weyl defines the kernel
$$
\psi_{\pm}^{\alpha}(t)=\sum_{k=-\infty,k\not=0}^{+\infty}\frac{e^{ikt}}{(\pm ik)^{\alpha}}=2\sum_{k=1}^{\infty}\frac{\cos(kt\mp\alpha\frac{\pi}{2})}{k^\alpha}.
$$
Thus the so called Marchaud-Weyl derivative is defined as
\begin{equation}\label{wewewe}
{\bold D}^{(\alpha)}_{\pm}f(x)=\frac{1}{2\pi}\int_0^{2\pi}\left(f(x)-f(x-t)\right)\frac{d}{dt}\psi^{1-\alpha}_\pm(t) dt.
\end{equation}
If the function $f$ is $2\pi$ periodic then it results that
$$
\frac{\alpha}{\Gamma(1-\alpha)}\int_0^{+\infty}\frac{f(x)-f(x-t)}{t^{1+\alpha}}dt=\frac{1}{2\pi}\int_0^{2\pi}\left(f(x)-f(x-t)\right)\frac{d}{dt}\psi^{1-\alpha}_\pm(t) dt.
$$
Previous relationship has to be correctly expressed in the following sense, see Lemma 19.4 in \cite{SKM}. 
\begin{prop}[Lemma 19.4, \cite{SKM}]\label{SamkoWeyl_Marchaud_coinc} 
For every $f\in L^p(0,2\pi),$ $1\geq p< +\infty$  the following limits converge for almost every $x\in (0,2\pi)$ simultanuously
$$
\lim_{\epsilon \to 0^+}\frac{1}{2\pi}\int_\epsilon^{2\pi}\left(f(x)-f(x-t)\right)\frac{d}{dt}\psi^{1-\alpha}_\pm(t) dt,
$$
$$
\lim_{\epsilon \to 0^+}\frac{\alpha}{\Gamma(1-\alpha)}\int_\epsilon^{+\infty}\frac{f(x)-f(x-t)}{t^{1+\alpha}}dt
$$
and
\begin{equation*}
\begin{split}
{\bold D}^{(\alpha)}_{+}f(x)&=\lim_{\epsilon \to 0^+}\frac{1}{2\pi}\int_\epsilon^{2\pi}\left(f(x)-f(x-t)\right)\frac{d}{dt}\psi^{1-\alpha}_\pm(t) dt\\
&=\lim_{\epsilon \to 0^+}\frac{\alpha}{\Gamma(1-\alpha)}\int_\epsilon^{+\infty}\frac{f(x)-f(x-t)}{t^{1+\alpha}}dt={\bold D}^{\alpha}_{+}f(x).
\end{split}
\end{equation*}
\end{prop}
By the way, concerning the parallel situation for the fractional Laplace operator on the torus, we point out \cite{Stinga_Roncal}, where similar results to the Proposition \ref{SamkoWeyl_Marchaud_coinc} have been proved.

We do not enter in the details concerning the question whether the definition of this type of fractional derivative has been invented by Weyl or Marchaud, or maybe by Riemann himself indeed as Weyl seems to suggest in the introduction of his paper \cite{Weyl}. Nevertheless the formula (\ref{Weyl1}) appeared in \cite{Weyl}, as already written, ten years before the Marchaud thesis. In any case, Marchaud correctly cites the Weyl's paper \cite{Weyl}. More precisely Marchaud at p. 50 of his thesis acknowledges to  Weyl  to have obtained the result in the case of dimension $n=1,$ by referring to the representation (\ref{Marchaud_der}). In addition Marchaud also admits that Weyl's approach was more powerful with respect to the one established by Montel, see \cite{Montel}. Montel approach used polynomial approximation, as Marchaud stated. On the contrary, Marchaud remarked, that Weyl's approach is more direct.  In \cite{SKM}, see , see XXXIII, the authors faced indirectly that question in the note dedicated to the historical outline of the subject. There, they explained that formula  (\ref{Marchaud_der}) appeared earlier in \cite{Weyl} by accident. Nevertheless, they concluded that Weyl did not develop his idea, as on the contrary Marchaud did in \cite{Marchaud_T} and \cite{Marchaud}. It would be interesting to know if any interaction between Weyl and Marchaud happened. In any case the Weyl's paper \cite{Weyl} is not one of the most cited among all the important results obtained by Weyl during his fruitful career.

\section{Basic ideas}\label{ideasMW} 
If we compare the Marchaud derivative  with respect to the Riemann-Liouville one, we immediately realize that in the latter one the classical derivative operator appears, while in the first one it does not.
 This is one of the key point that Marchaud's  definition put in evidence. That is, Marchaud derivative avoids to apply the classical derivative after an integration in order to define the fractional operator. In a sense, this approach recalls the one that has characterized the Sobolev's approach, see for instance \cite{Sobolev}  and, in 
in a sense, it could be considered as precursive of the notion of weak solution to a PDEs. In fact, roughly saying, we recall that  Sobolev's approach is based on the  integration of the both sides of an equation. In this way we reduced to look for functions that satisfy the obtained integral equation.
 
In this order of ideas, in Riemann-Liouville definition, it still appears the classical derivation. On the contrary, in Marchaud derivative we simply recognize a singular integral where the reminiscence of the derivative is given by the kernel that multiplies the difference between the values of the function in two points. On the other hand, Marchaud's definition includes the Riemann-Liouville's one when the initial point is $-\infty$ and the functions are sufficiently smooth. We come back to this aspect later on in the section. From a philosophical point of view, Marchaud derivative seems to put in evidence its non-local character. On the contrary, in the initial historical approach described by Riemann-Liouville derivative, the classical derivative operator, that is a local object, still remains. 

For instance, by considering a function defined in all of $\mathbb{R}$ and having a minimal smoothness we, in principle, can modify its definition locally, for example simply changing its derivative in a {\it small} set of points. Nevertheless, the remaining part of the function is not affected by this modification. On the contrary the Marchaud derivative, but also the Riemann-Liouville's one, even if in a spurious way, determines a quantity that heavily depends on the modified function. This fact is evident thanks to the presence of the integral operator. Summarizing, by modifying the given function even only in a small set, the value of the fractional derivative will change, in general, in all the points where this fractional derivative will be evaluated.  

Now we comment separately Marchaud derivative and Weyl derivative.

\subsection{Marchaud derivative}
Marchaud derivative acts like an operator that associates to a function a new function that in general does not maintain local properties like the differential (of the function) do far away to the set where the function has been modified. Nevertheless, this operator, the fractional one, in a sense, still contains the classical derivative. Indeed the classical derivative materializes as a particular (let say like an exception) case who realizes when the order of the fractional derivative goes to an integer. This focusing phenomenon is particularly interesting.

In order to clarify this remark let us consider the definition (\ref{exampleMarchaud}), in the case $p=2$ and $\alpha=1.$ Then we obtain:
\begin{equation}\label{Marchaudf}
\begin{split}
{\bold D}^{\alpha}f(x)\int_{0}^{\infty}s^{-2}(1-e^{-s})^2ds=\int_0^\infty s^{-\alpha-1}\left(f(x)-f(x-s)+f(x-2s)\right)ds,
 \end{split}
\end{equation}
and since
\begin{equation*}
\begin{split}
&\int_{0}^{\infty}s^{-2}(1-e^{-s})^2ds
=2\log 2,
\end{split}
\end{equation*}
we get:
\begin{equation*}
\begin{split}
{\bold D}f(x)=\frac{1}{2\log 2}\int_0^\infty \frac{f(x)-2f(x-s)+f(x-2s)}{s^{2}}ds.
 \end{split}
\end{equation*}
It is worth to say that here we have the value of the classical derivative in a point represented via an integral! Let say: from the global to the local. How to explain this fact?
We remark  that if $f$ is a $C^2$ function with compact support or even $f\in \mathcal{S}(\mathbb{R})$, then  
\begin{equation*}
\begin{split}
&\mathcal{F}\left(\int_0^\infty \frac{f(x)-2f(x-s)+f(x-2s)}{s^{2}}ds\right)=\mathcal{F}f(\xi)\int_0^{\infty}\frac{(1-e^{-is\xi})^2}{s^{2}}ds\\
&=2\log2(i\xi)\mathcal{F}f(\xi)=2\log2\mathcal{F}(f').
 \end{split}
\end{equation*}
This implies, recalling that the Fourier transform is invertible on Schwartz space $\mathcal{S}(\mathbb{R}),$  that for every $f\in \mathcal{S}(\mathbb{R})$ formula (\ref{Marchaudf}) truly gives a representation of the derivative of a function in a point. We shall come back in Section  \ref{general_setting} on this fact. On the other hand, the relationship (\ref{Marchaudf}) is correctly defined in a larger space of functions with respect to $\mathcal{S}(\mathbb{R}).$

In the case $p=1,$ $\alpha<1$
\begin{equation*}
\begin{split}
{\bold D}^{\alpha}f(x)\int_{0}^{\infty}s^{-1-\alpha}(1-e^{-s})ds=\int_0^\infty s^{-\alpha-1}\left(f(x)-f(x-s)\right)ds,
 \end{split}
\end{equation*}
but 
$$
\int_{0}^{\infty}s^{-1-\alpha}(1-e^{-s})ds=[-\alpha^{-1} s^{-\alpha}(1-e^{-s})]_{s=0}^{s=\infty}+\alpha^{-1}\int_{0}^{\infty}s^{-\alpha}e^{-s}ds=\alpha^{-1}\Gamma(1-\alpha).
$$
As a consequence
\begin{equation}\label{Marchaud_der}
{\bold D}^{\alpha}f(x)=\frac{\alpha}{\Gamma(1-\alpha)}\int_0^\infty \frac{f(x)-f(x-s)}{s^{1+\alpha}}ds
\end{equation}
and even in this case, the easier case among Marchaud derivatives concerning  the function $f$ for $\alpha\in (0,1),$ we can read the non-locality of this definition and, in addition, for sufficiently smooth functions,
$$
\lim_{\alpha\to1^-}{\bold D}^{\alpha}f(x)=Df(x).
$$
In this case an important role is played by the normalizing constant $\frac{\alpha}{\Gamma(1-\alpha)}$ that multiplies the integral in the definition of Marchaud derivative.

The fact that, for sufficiently "good" functions, the fractional derivative ${\bold D}^{\alpha}f$ coincides with the Riemann-Liouville derivative
$$
\mathcal{D}^\alpha f(x)=\frac{1}{\Gamma(1-\alpha)}\frac{d}{dx}\int_{-\infty}^x\frac{f(t)}{(x-t)^{\alpha}},
$$
can be checked straightforwardly.
Moreover, the definition given by Marchaud can be applied even for functions that may growth at infinity less than $\alpha.$ On the contrary the definition of Liouville derivative is less flexible since it does not admit (see p.XXXIII \cite{SKM}) to be applied to constant functions. 
 
 Let us check that Marchaud derivative ${\bold D}^{\alpha}_+f$
coincides with Riemann-Louville derivative from the right. In fact, since
$$
\mathcal{D}^\alpha_+f(x)=\frac{1}{\Gamma(1-\alpha)}\frac{d}{dx}\int_0^{+\infty}\frac{f(x-t)}{t^\alpha}dt,
$$ 
and supposing that $f\in C^1(\mathbb{R})$ and $f=o(|x|^{\alpha-1-\epsilon}),$  $x\to +\infty$ for $\epsilon>0,$
then by Lebesgue dominated convergence theorem first and then integrating by parts, we get:
\begin{equation}
\begin{split}
\mathcal{D}^\alpha_+f(x)&=\frac{1}{\Gamma(1-\alpha)}\int_0^{+\infty}\frac{f'(x-t)}{t^\alpha}dt=\frac{\alpha}{\Gamma(1-\alpha)}\int_0^{+\infty}f'(x-t)\left(\int_{t}^{\infty}\tau^{-\alpha-1}d\tau\right)dt\\
&=\frac{\alpha}{\Gamma(1-\alpha)}\lim_{\epsilon\to 0^+}\left\{[-f(x-t)\left(\int_{t}^{+\infty}\tau^{-\alpha-1}d\tau\right)]_{t=\epsilon}^{+\infty}-\int_\epsilon^{+\infty}\frac{f(x-t)}{t^{1+\alpha}}dt\right\}\\
&=\frac{\alpha}{\Gamma(1-\alpha)}\lim_{\epsilon\to 0^+}\left\{f(x-\epsilon)\left(\int_{\epsilon}^{+\infty}\tau^{-\alpha-1}d\tau\right)-\int_\epsilon^{+\infty}\frac{f(x-t)}{t^{1+\alpha}}dt\right\}\\
&=\frac{\alpha}{\Gamma(1-\alpha)}\lim_{\epsilon\to 0^+}\left\{\int_\epsilon^{+\infty}\frac{f(x)-f(x-t)}{t^{1+\alpha}}dt+(f(x-\epsilon)-f(x))\int_{\epsilon}^{+\infty}\tau^{-\alpha-1}d\tau\right\}=\\
&=\frac{\alpha}{\Gamma(1-\alpha)}\int_0^{+\infty}\frac{f(x)-f(x-t)}{t^{1+\alpha}}dt={\bold D}^{\alpha}_+f(x),
\end{split}
\end{equation}
because  there exists $\eta\in ]x-\epsilon, x[$ such that, as $\epsilon \to 0:$
$$
|\alpha(f(x-\epsilon)-f(x))\int_{\epsilon}^{+\infty}\tau^{-\alpha-1}d\tau|=|\alpha f'(x-\eta)|\leq \sup_{\tau\in[x-\epsilon, x]}|f'(\tau)| \epsilon^{1-\alpha}\to 0.
$$
Thus, from this point of view, Marchaud derivative is a sort of weaker version of the Riemann-Liouville derivative. 

For example, constants satisfy ${\bold D}^{\alpha}_+f(x)=0$ in the Marchaud sense, 
even if we can not consider, in all of $\mathbb{R},$  the Riemann-Liouville derivative of a constant. In fact the parallel integral is divergent. This is, of course, absolutely unpleasant! Indeed, both Marchaud and Weyl were motivated also from this fact in order for looking for a different type of definition of fractional derivative. 

We also think that Marchaud derivative as some further properties that have to be better understood in its application. In order to focus one of this aspect, we remark, see also \cite{Fer_H}, that the sum of the two Marchaud derivatives (${\bold D}^{\alpha}_+f$ and ${\bold D}^{\alpha}_-f$) gives, in a sense, the Riesz derivative in one dimension, namely the fractional Laplace operator in dimension $1.$ More precisely:
\begin{equation*}
\begin{split}
{\bold D}^{\alpha}_+f(x)+{\bold D}^{\alpha}_-f(x)=\frac{\alpha}{\Gamma(1-\alpha)}\int_0^{+\infty}\frac{2f(x)-f(x-t)-f(x+t)}{t^{1+\alpha}}dt
\end{split}
\end{equation*}
or
\begin{equation*}
\begin{split}
{\bold D}^{\alpha}_+f(x)+{\bold D}^{\alpha}_-f(x)&=\frac{\alpha}{\Gamma(1-\alpha)}\left(\int_0^{+\infty}\frac{f(x)-f(x-t)}{t^{1+\alpha}}dt+\int_0^{+\infty}\frac{f(x)-f(x+t)}{t^{1+\alpha}}dt\right)\\
&=\frac{\alpha}{\Gamma(1-\alpha)}\left(\int_0^{+\infty}\frac{f(x)-f(x-t)}{t^{1+\alpha}}dt+\int^0_{-\infty}\frac{f(x)-f(x-t)}{|t|^{1+\alpha}}dt\right)\\
&=\frac{\alpha}{\Gamma(1-\alpha)}\int_{-\infty}^{+\infty}\frac{f(x)-f(x-t)}{|t|^{1+\alpha}}dt=\frac{\alpha}{\Gamma(1-\alpha)}\int_{-\infty}^{+\infty}\frac{f(x)-f(\xi)}{|x-\xi|^{1+\alpha}}d\xi\\
&=\frac{\alpha}{c(1,\frac{\alpha}{2})\Gamma(1-\alpha)}(-\frac{d^2}{dx^2})^{\frac{\alpha}{2}}f(x),
\end{split}
\end{equation*}
where $c(1,\frac{\alpha}{2})$ is the normalizing constant associated with the fractional Laplace operator $(-\frac{d^2}{dx^2})^{\frac{\alpha}{2}},$ and whose value we will recall later on in this paper. This fact was implicitly remarked in \cite{SKM} and it seems that it can be connected with the different type of variable considered. In case of ${\bold D}^{\alpha}_+f$ and ${\bold D}^{\alpha}_-f$ the only one variable in $\mathbb{R}$ has a privileged direction in the two definitions of fractional derivative. For instance, we can think to it as the time variable. On the contrary considering the fractional Laplace operator $$(-\frac{d^2}{dx^2})^{\frac{\alpha}{2}}f$$ in $\mathbb{R}$ (the same also in $\mathbb{R}^n$) there is not any privileged direction. Namely the space (in this case $\mathbb{R}$) is homogeneous. So that previous connection is particularly interesting.

\subsection{Weyl derivative}
As far as Weyl's approach concerned, the relationship between spectral theory and fractional derivative is explicit. Indeed,  supposing of working with a $2\pi$-periodic function having zero average, it is well known that the associated Fourier series is
$$
\sum_{k=-\infty}^{+\infty}c_ke^{ikx},
$$
where of course $\{c_k\}_{k\in \mathbb{Z}}$ denotes the sequence of Fourier coefficients.

Then, by computing formally the derivative of this series, we obtain
$$
\sum_{k=-\infty}^{+\infty}c_k(ik)e^{ikx}.
$$
It is obvious that defining a new function for a fixed $\alpha<1,$ as
$$
\sum_{k=-\infty}^{+\infty}\frac{c_k}{(ik)^\alpha}e^{ikx},
$$
we formally obtain taking then a derivative we obtain
\begin{equation}\label{wewe}
D\left(\sum_{k=-\infty}^{+\infty}\frac{c_k}{(ik)^\alpha}e^{ikx}\right)=\sum_{k=-\infty}^{+\infty}\frac{c_k}{(ik)^{\alpha-1}}e^{ikx}.
\end{equation}
In this way, Weyl defines the parallel fractional integral so that it is natural to define the fractional derivative of $f$ as
$$
\sum_{k=-\infty, }^{+\infty}c_k(ik)^\alpha e^{ikx}.
$$
On the other hand, we recall that given two periodic functions $f,g$ the new function
$$
\frac{1}{2\pi}\int_0^{2\pi}g(t)f(x-t)dt
$$
is represented by the Fourier series
$$
\sum_{k=-\infty}^{\infty}g_kc_ke^{ikx},
$$
where $\{g_k\}_{k\in\mathbb{Z}}$ and $\{c_k\}_{k\in\mathbb{Z}}$ are the respective  Fourier coefficients.

As a consequence, considering 
$$
\sum_{k=-\infty}^{+\infty}\frac{c_k}{(ik)^\alpha}e^{ikx},
$$
as representing the Fourier series of an integral like the following one,
$$
\frac{1}{2\pi}\int_0^{2\pi}g(t)f(x-t)dt,
$$
we desume that previous integral has to be written in the following form:

$$
\frac{1}{2\pi}\int_0^{2\pi}f(x-t)\left(\sum_{k=-\infty,k\not=0}^{+\infty}\frac{e^{ikt}}{(ik)^\alpha}\right)dt.
$$
Since it can prove that,  see \cite{SKM}, that
$$
\sum_{k=-\infty,k\not=0}^{+\infty}\frac{e^{ikt}}{(ik)^\alpha}=2\sum_{k=1}^{\infty}\frac{\cos(kt-\alpha\frac{\pi}{2})}{k^\alpha}.
$$
Then, denoting the kernel
$$
\psi_{+}^{\alpha}(t):=\sum_{k=-\infty,k\not=0}^{+\infty}\frac{e^{ikt}}{(ik)^\alpha},
$$
Weyl obtains the fractional integral
$$
I^{(\alpha)}_+f(x)=\frac{1}{2\pi}\int_0^{2\pi}f(x-t)\psi_{+}^{\alpha}(t)dt.
$$ 
At this point, by recalling (\ref{wewe}), Weyl defines the fractional derivative as
$$
\mathcal{D}_+^{(\alpha)}(x)=D\left(I^{(1-\alpha)}_+f\right)(x).
$$
This definition corresponds to the Weyl-Riemann-Lioville version of this derivative, see \cite{SKM} for the details. Then taking formally the derivative Weyl obtains
the Weyl-Marchaud derivative, see  also(\ref{wewewe}), discussed in Section \ref{Weyl_approach}:
$$
{\bold D}^{(\alpha)}_{+}f(x)=\frac{1}{2\pi}\int_0^{2\pi}\left(f(x)-f(x-t)\right)\frac{d}{dt}\psi^{1-\alpha}_+(t) dt.
$$
Of course the case concerning ${\bold D}^{(\alpha)}_{-}f$ is analogous to the one just described for ${\bold D}^{(\alpha)}_{+}f$.


\section{General setting of Marchaud derivative and some further remarks} \label{general_setting}

The definition of Marchaud derivative, as it is known since \cite{SKM}, can be extended to all $\alpha>0$ in the following way, see \cite{SKM} and \cite{Samko}. Let $l\in \mathbb{N},$  $l\geq 1$ and $\alpha<l.$ We define for every $f\in \mathcal{S}(\mathbb{R})$
$$
\bold{D}^\alpha_{\pm}f(x)=\frac{1}{\chi(\alpha,l)}\int_{0}^{+\infty}\frac{\Delta^l_{\pm \tau}f(x)}{\tau^{1+\alpha}}d\tau,
$$
where
$$
\chi(\alpha,l)=\Gamma(-\alpha)A_l(\alpha)=\int_0^{+\infty}\frac{(1-e^{-t})^l}{t^{1+\alpha}}dt,
$$
\begin{equation}\label{A_coeffi}
A_l(\alpha)=\sum_{k=0}^l(-1)^k\binom{l}{k}k^\alpha,
\end{equation}
and
$$
\Delta^l_{\pm \tau}f(x)=\sum_{k=0}^l(-1)^k\binom{l}{k}f(x\mp k\tau).
$$
Of course this definition can be generalized to the case of functions with several variables having a nice behavior both at infinity and locally, simply considering for every $\xi\in \mathbb{R}^n,$ $\xi\not=0$  and for every $\alpha>0$ and $l\in \mathbb{N},$  $l\geq 1$ such that $\alpha<l$ and defining:
$$
\mathcal{D}^{\alpha)}_{\pm, \xi}f(x)=\frac{1}{\chi(\alpha,l)}\int_{0}^{+\infty}\frac{\Delta^l_{\pm\tau,\xi}f(x)}{\tau^{1+\alpha}}d\tau,
$$
where 
$$
\Delta^l_{\pm\tau,\xi}f(x)=\sum_{k=0}^l(-1)^k\binom{l}{k}f(x\mp k\tau\xi).
$$
It is worth to say that $\lim_{\alpha\to l^-}\mathcal{D}^{\alpha)}_{\pm, \xi}f(x)=\pm {D}^{l}_{\xi}f(x),$
in the local (classical sense) and  $\lim_{\alpha\to (l-1)^+}\mathcal{D}^{\alpha)}_{\pm, \xi}f(x)=\pm D^{l-1}_{\xi}f(x),$
where $\mathcal{D}^{0)}_{\xi}=I,$  whenever $f$ is sufficiently smooth (for example in $\mathcal{S}(\mathbb{R}^n)$).

In this way it is possible to consider interesting representation of local operators. For example denoting by $e_i$ the vector of the canonic base of $\mathbb{R}^n,$ for every $i=1,\dots,n$ we get
\begin{equation}\begin{split}
\frac{\partial f(x)}{\partial x_i}=\frac{1}{\chi(1,2)}\int_{0}^{+\infty}\frac{\Delta^2_{\tau,e_i}f(x)}{\tau^{2}}d\tau.
\end{split}
\end{equation}
and
\begin{equation}\begin{split}
\frac{\partial^2 f(x)}{\partial x_i^2}=\frac{1}{\chi(2,3)}\int_{0}^{+\infty}\frac{\Delta^3_{\tau,e_i}f(x)}{\tau^{3}}d\tau.
\end{split}
\end{equation}
As a consequence for every $f\in \mathcal{S}(\mathbb{R}^n):$
\begin{equation}\begin{split}
\sum_{i=1}^n\frac{\partial}{\partial x_i}\frac{\partial f(x)}{\partial x_i}=\Delta f(x)=\frac{1}{\chi(1,2)}\int_{0}^{+\infty}\frac{\sum_{i=1}^n\left(\frac{\partial f(x)}{\partial x_i }+\sum_{k=1}^2(-1)^{k+1}k\binom{2}{k}\tau\frac{\partial f(x-k\tau\xi)}{\partial x_i}\right)}{\tau^{2}}d\tau.
\end{split}
\end{equation}
and
\begin{equation}\begin{split}
\Delta f(x)=\sum_{i=1}^n\frac{\partial^2 f(x)}{\partial x_i^2}=\frac{1}{\chi(2,3)}\int_{0}^{+\infty}\frac{\sum_{i=1}^n\Delta^3_{\tau,e_i}f(x)}{\tau^{3}}d\tau,
\end{split}
\end{equation}
that is
\begin{equation}\begin{split}
\Delta f(x)=\frac{1}{\chi(2,3)}\int_{0}^{+\infty}\frac{\sum_{k=0}^3(-1)^k\binom{3}{k}\sum_{i=1}^nf(x-k\tau e_i)}{\tau^{3}}d\tau,
\end{split}
\end{equation}
From Liouville Theorem it is well known that there exists a unique function $f\in\mathcal{S}(\mathbb{R}^n)$ such that $\Delta f(x)=0$ in $\mathbb{R}^n$ that is $f=0.$ 
Thus the unique function $f\in\mathcal{S}(\mathbb{R}^n)$ that satisfies 
$$
\frac{1}{\chi(2,3)}\int_{0}^{+\infty}\frac{\sum_{i=1}^n\Delta^3_{\tau,e_i}f(x)}{\tau^{3}}d\tau=0
$$
 has to be  $f=0.$
 
 About the properties of the Marchaud derivative we like to remind that for every function $f\in\mathcal{S}(\mathbb{R})$
 \begin{equation*}
\begin{split}
{\bold D}f(x)=\frac{1}{2\log 2}\int_0^\infty \frac{f(x)-2f(x-s)+f(x-2s)}{s^{2}}ds.
 \end{split}
\end{equation*}
On the other hand for every $f,g\in \mathcal{S}(\mathbb{R}),$ $fg\in\mathcal{S}(\mathbb{R}),$ so that
 \begin{equation*}
\begin{split}
&{\bold D}(fg)(x)=\frac{1}{2\log 2}\int_0^\infty \frac{f(x)g(x)-2f(x-s)g(x-s)+f(x-2s)g(x-2s)}{s^{2}}ds
 \end{split}
\end{equation*}
on the other hand we know that
${\bold D}(fg)(x)={\bold D}(f)(x)g(x)+{\bold D}(g)(x)f(x).$ Then as a byproduct we obtain the following formula for every $f,g\in \mathcal{S}(\mathbb{R})$
\begin{equation*}
\begin{split}
&\int_0^\infty \frac{f(x)g(x)-2f(x-s)g(x-s)+f(x-2s)g(x-2s)}{s^{2}}ds\\
&=\int_0^\infty \frac{f(x)-2f(x-s)+f(x-2s)}{s^{2}}dsg(x)+\int_0^\infty \frac{g(x)-2g(x-s)+g(x-2s)}{s^{2}}dsf(x).\\
 \end{split}
\end{equation*}
Nevertheless, for instance, for every $\alpha\in (0,1)$ and for every $f,g\in \mathcal{S}(\mathbb{R})$ we get that $fg\in\mathcal{S}(\mathbb{R})$ and
 \begin{equation*}
\begin{split}
&{\bold D}^\alpha (fg)(x)=\frac{\alpha}{\Gamma(1-\alpha)}\int_0^{\infty}\frac{f(x)g(x)-f(x-t)g(x-t)}{t^{1+\alpha}}dt\\
&=\frac{\alpha}{\Gamma(1-\alpha)}\int_0^\infty\frac{f(x)-f(x-t)}{t^{1+\alpha}}dtg(x)+\frac{\alpha}{\Gamma(1-\alpha)}\int_0^\infty \frac{f(x-t)(g(x)-g(x-t))}{t^{1+\alpha}}dt\\
&={\bold D}^\alpha f(x)g(x)+{\bold D}^\alpha g(x)f(x)-\frac{\alpha}{\Gamma(1-\alpha)}\int_0^\infty \frac{(f(x)-f(x-t))(g(x)-g(x-t))}{t^{1+\alpha}}dt.
 \end{split}
\end{equation*}
This remark implies that  the usual differential rule for the product of two functions does not hold. Nevertheless
 \begin{equation*}
\begin{split}
\frac{\alpha}{\Gamma(1-\alpha)}\int_0^\infty \frac{(f(x)-f(x-t))(g(x)-g(x-t))}{t^{1+\alpha}}dt\to 0
 \end{split}
\end{equation*}
whenever $\alpha\to 1^-.$ In fact
 \begin{equation*}
\begin{split}
&\frac{\alpha}{\Gamma(1-\alpha)}\int_0^\infty \frac{(f(x)-f(x-t))(g(x)-g(x-t))}{t^{1+\alpha}}dt\\
&=\frac{\alpha}{\Gamma(1-\alpha)}\int_0^\eta \frac{(f(x)-f(x-t))(g(x)-g(x-t))}{t^{1+\alpha}}dt\\
&+\frac{\alpha}{\Gamma(1-\alpha)}\int_\eta^\infty \frac{(f(x)-f(x-t))(g(x)-g(x-t))}{t^{1+\alpha}}dt\\
&=\frac{\alpha}{\Gamma(1-\alpha)}Df(x)Dg(x)\int_0^\eta t^{1-\alpha}dt+o(\eta^{2-\alpha})\\
&+\frac{\alpha}{\Gamma(1-\alpha)}\int_\eta^\infty \frac{(f(x)-f(x-t))(g(x)-g(x-t))}{t^{1+\alpha}}dt\to 0,
 \end{split}
\end{equation*}
whenever $\alpha\to 1^-,$ because $\frac{\alpha}{\Gamma(1-\alpha)}\to 0,$ and there exists a positive constant such that
$$
\left|\int_\eta^\infty \frac{(f(x)-f(x-t))(g(x)-g(x-t))}{t^{1+\alpha}}dt\right|\leq M
$$
uniformly for every $\alpha\in (0,1)$ and for every fixed $\eta>0.$

 In this way we obtain one more time the classical rule for the usual derivative of order one, because ${\bold D}^\alpha (fg)(x)\to {\bold D}(fg)(x)$  and ${\bold D}^\alpha f(x)g(x)+{\bold D}^\alpha g(x)f(x)\to{\bold D}f(x)g(x)+{\bold D}g(x)f(x)$ if $\alpha\to 1^-.$
 
 This behavior is heuristically clear thinking to the fractional operator as a nonlocal object. That is, the fractional derivative in a point {\it measure} something that depends on all the values of the function before that point. Thus it is in a sense expected that for this type of operator a term depending on the interplay of the quantity associated with the fractional derivative of the functions acting, has to appear. In the special case of $\alpha\to 1^-$ this third term appears with value $0$ thanks to the locality of the quantity expressed by the classical derivative of order one. The Marchaud derivative of order $\alpha,$ it rescales with the law $\lambda^\alpha.$ In fact we have that for every $f\in \mathcal{S}(\mathbb{R}),$ the function $x\to f(\lambda x)=f_\lambda(x)$ has the following behavior with respect to the Marchaud fractional derivative
  \begin{equation*}
\begin{split}
{\bold D}^\alpha f_\lambda(x)&= \frac{\alpha}{\Gamma(1-\alpha)}\int_0^{\infty}\frac{f(\lambda x)-f(\lambda(x-t))}{t^{1+\alpha}}dt\\
&= \frac{\alpha}{\Gamma(1-\alpha)}\int_0^{\infty}\frac{f(\lambda x)-f(\lambda x-\lambda t))}{t^{1+\alpha}}dt=\lambda^\alpha\frac{\alpha}{\Gamma(1-\alpha)}\int_0^{\infty}\frac{f(\lambda x)-f(\lambda x-\tau)}{\tau^{1+\alpha}}d\tau\\
&=\lambda^\alpha{\bold D}^\alpha f(\lambda x).
 \end{split}
\end{equation*}
Remark that also with respect to a different representation of the Marchaud fractional derivative, let see the case $\alpha<2,$ we get the same rescaling law
  \begin{equation*}
\begin{split}
{\bold D}^\alpha f_\lambda(x)&= \frac{1}{\chi(\alpha,2)}\int_0^{\infty}\frac{f(\lambda x)-2f(\lambda(x-t))+f(\lambda(x-2t))}{t^{1+\alpha}}dt\\
&= \frac{1}{\chi(\alpha,2)}\int_0^{\infty}\frac{f(\lambda x)-2f(\lambda x-\lambda t))+f(\lambda x-2\lambda t))}{t^{1+\alpha}}dt\\
&=\lambda^\alpha\frac{1}{\chi(\alpha,2)}\int_0^{\infty}\frac{f(\lambda x)-2f(\lambda x-\tau)+f(\lambda x-2\tau)}{\tau^{1+\alpha}}d\tau\\
&=\lambda^\alpha{\bold D}^\alpha f(\lambda x).
 \end{split}
\end{equation*}

The definition of Marchaud derivative makes sense for a larger class of functions, with respect to the set $\mathcal{S}(\mathbb{R}).$ For instance, all the constants have Marchaud derivative zero. The exponential function $x\to e^{\lambda x}$ does not belong to 
 $\mathcal{S}(\mathbb{R}).$ Nevertheless, for $\lambda\geq 0$ (here we are using the Marchaud derivative ${\bold D}^\alpha_+$ but for the sake of simplicity we omit of writing the sign $+$) the function $e^{\lambda x}$ has Marchaud derivative and
 \begin{equation*}
\begin{split}
 {\bold D}^\alpha e^{\lambda x}&=\frac{\alpha}{\Gamma(1-\alpha)}\int_0^{\infty}\frac{e^{\lambda x}-e^{\lambda(x-t)}}{t^{1+\alpha}}dt=e^{\lambda x}\frac{\alpha}{\Gamma(1-\alpha)}\int_0^{\infty}\frac{1-e^{-\lambda t}}{t^{1+\alpha}}dt\\
 &=\lambda^{\alpha}e^{\lambda x}\frac{\alpha}{\Gamma(1-\alpha)}\int_0^{\infty}\frac{1-e^{-\tau}}{\tau^{1+\alpha}}d\tau=\lambda^{\alpha}e^{\lambda x},
 \end{split}
\end{equation*}
because, integrating by parts we obtain:
 \begin{equation*}
\begin{split}
\int_0^{\infty}\frac{1-e^{-\tau}}{\tau^{1+\alpha}}d\tau=\frac{1}{\alpha}\int_0^{\infty}\frac{e^{-\tau}}{\tau^{\alpha}}d\tau=\frac{\Gamma(1-\alpha)}{\alpha}.
 \end{split}
\end{equation*}
As a consequence $e^{\lambda x}$ is solution of the fractional differential equation
$$
{\bold D}^\alpha f(x)=\lambda^\alpha f(x).
$$

\section{Fractional Laplace operator }\label{fractional_recall}
The fractional Laplace operator can be represented in several way. We should have to cite the  contribution of many authors, we recall for instance \cite{Riesz1},  \cite{Riesz2}, \cite{Riesz3} \cite{Riesz4} \cite{Landkof}, \cite{Stein}. Using the Fourier transform, for every $s\in (0,1)$ and for every $u\in \mathcal{S}(\mathbb{R}^n)$ the fractional Laplace operator is usually defined as
$$
(-\Delta)^su=\mathcal{F}^{-1}(\|\xi\|^{2s}\mathcal{F})u.
$$
As a consequence for every $u\in L^2(\mathbb{R}^n)$ if $\|\xi\|^{2s}\mathcal{F}u\in L^2(\mathbb{R}^n),$ then the fractional Laplace operator is defined by  $\mathcal{F}^{-1}(\|\xi\|^{2s}\mathcal{F})u.$

On the other hand  for every $u\in \mathcal{S}(\mathbb{R}^n)$ and $s\in (0,1)$ we can define the operator
$$
\mathcal{L}_s u(x)=c(\alpha,n)\int_{\mathbb{R}^n}\frac{f(x)-f(y)}{\|x-y\|^{n+2s}}dy:=\lim_{\epsilon\to 0}c(s,n)\int_{\mathbb{R}^n\setminus B_\epsilon(x)}\frac{f(x)-f(y)}{\|x-y\|^{n+2s}}dy
$$
where $c(\alpha,n)$ is a normalizing constant, then $\mathcal{L}_s=(-\Delta)^s$ and 
$$
c(s,n)=\left(\int_{\mathbb{R}^n}\frac{1-\cos(\xi_1)}{\|\xi\|^{n+2s}}d\xi\right)^{-1}.
$$
In addition, see \cite{DPV}, if $n>1,$ we get:
$$
\lim_{s\to{1^-}}\frac{\omega_{n-1}c(s,n)}{4ns(1-s)}=1
$$
and
$$
\lim_{s\to{0^+}}\frac{\omega_{n-1}c(s,n)}{2s(1-s)}=1.
$$
In addition in Lemma 5, \cite{StiTor}, previous constant has been surprisingly computed in a precise way so that it results:
$$
c_{s,n}=\frac{4^s\Gamma(\frac{n}{2}+s)}{\pi^{\frac{n}{2}}\Gamma(-\sigma)}.
$$
We recall that in \cite{SKM} a different expression of the fractional Laplace operator has been given, introducing a different constant of normalization and considering a more general situation. In fact for every $f\in \mathcal{S}(\mathbb{R}^n)$ and $\alpha>0,$  $l\in \mathbb{N},$ $n\geq 1,$ $\alpha<l$ we
may define the following operator:
$$
(-S)^{\frac{\alpha}{2}}f(x)=\frac{\sin(\alpha\frac{\pi}{2})}{\beta_n(\alpha)A_l(\alpha)}\int_{\mathbb{R}^n}\frac{\Delta^l_yf(x)}{\|y\|^{n+\alpha}}dy,
$$
where $A_l(\alpha)$ is defined in (\ref{A_coeffi}),
$$
\beta_n(\alpha)=\frac{\pi^{1+\frac{n}{2}}}{2^\alpha\Gamma(1+\frac{\alpha}{2})\Gamma(\frac{n+\alpha}{2})}
$$
and
$$
\Delta^l_yf(x)=\sum_{k=0}^l(-1)^k\binom{l}{k}f(x-ky).
$$
denotes the non-centered differences. Then in \cite{SKM}, see Lemma 25.3, it is possible to find the proof that $(-S)^{\alpha/2}=(-\Delta)^{\alpha/2}$ in $\mathcal{S}(\mathbb{R}^n).$

Another way of introducing the fractional Laplace operator can be done considering if $U:\mathbb{R}^n\times]0,+\infty[\to\mathbb{R}$  solution of the following nonlocal problem, 
\begin{equation*}
\left\{
\begin{array}{l}
\mbox{div}_{(x,y)}(y^{1-2s}\nabla U(x,y))=0,\quad\mbox{in}\quad \mathbb{R}^n\times]0,+\infty[\\
U(\cdot,0)=u,\quad x\in\mathbb{R}^n.
\end{array}
\right.
\end{equation*}
Then defining
$$
\mathcal{N}_su:=\lim_{y\to 0}y^{1-2s}\frac{\partial U(\cdot,y)}{\partial y}
$$
it results, possibly up to a multiplicative factor depending only on $s$ and $n$ to $\mathcal{N}_s,$ that  $(-\Delta)^s=(-S)^s=\mathcal{L}_s=\mathcal{N}_s$ for every $u\in \mathcal{S}(\mathbb{R}^n).$ Among the application of this extension approach we have the application to Carnot groups, see \cite{FerFra} and \cite{StiTor}.

In the next Section \ref{FerBucSec}, we shall discuss the relationship of the Marchaud derivative with respect to the previous representation of the fractional Laplace operator.  We recall however that, for the sake of completeness, the fractional Laplace operator may be represented also defining the operator
$$
\mathcal{A}_s=-\frac{s}{\Gamma(1-s)}\int_{0}^{+\infty}(e^{t\Delta}-Id)\frac{dt}{t^{1+s}},
$$
where $e^{t\Delta}$ denotes the heat semigroup generated by the Laplace operator $\Delta$ and it is also well known that defining the operator
$$
\mathcal{B}_s=c(s,n)\int_0^{+\infty}\lambda^sdE(\lambda),
$$
where $\{E(\lambda)\}_{\lambda\in [0,\infty[}$ is, as usual, the family of spectral projectors of the Laplace operator, we can conclude that, at least in $\mathcal{S}(\mathbb{R}^n),$ $(-\Delta)^s=(-S)^{s}=\mathcal{L}_s=\mathcal{N}_s=\mathcal{A}_s=\mathcal{B}_s.$
We conclude this section recalling \cite{BalaAV}, where the semigroup method has been introduced.  In \cite{StiTor} this approach has been developed and then generalized in \cite{GMS} to a very large class of operators.
The fractional Laplace operator in its  representation via an extension has been applied in \cite{CSS} for facing the regularity of the thin obstacle, see also \cite{SaSa}. In particular, we point out that this approach opened the way to a large number of papers in which  this idea applied to many other problems. Other applications of the fractional calculus to the geometric measure theory can be found, for instance, in \cite{CaVa}, \cite{AmDeMa}, \cite{CaSaVa} and  also coming to a very recently result \cite{CiFraGo} and \cite{FMPPS}, where the definition of nonlocal (fractional) perimeters is discussed. 
 For further insights to the properties of the fractional Laplace operator, in addition to \cite{Landkof}, \cite{SKM}, \cite{DPV}, and \cite{BiRaRa}, \cite{BuVa}, we point out also the very recent preprint \cite{Ga}.
\section{Extension approach for Marchaud derivative}\label{FerBucSec}

We described  here the simplest case given by $s={1}/{2}$ as follows. Let $\varphi:\mathbb{R}\to \mathbb{R}$ be a function in $\mathcal{S}(\mathbb{R})$ and $U$ be a solution of the problem
  \begin{equation} \label{pro1}\displaystyle
 \left\{ 
  \begin{split}
 &\frac{\partial U}{\partial t} =\frac{\partial^2 U}{\partial x^2}, & & (x,t) \in (0,\infty)\times\mathbb{R}\\
 		&U(0,t)=\varphi(t), & &t \in \mathbb{R}. 
\end{split}
		\right.
		\end{equation}
It is worth to remark that this is not the usual Cauchy problem for the heat operator. It is a heat conduction problem. 

Without extra assumptions, we can not expect to have a unique solution of the problem (\ref{pro1}), see \cite{tichonov}, Chapter 3.3. Anyhow, if we denote by $T_{1/2}$ the operator that associates to $\varphi$ the partial derivative ${\partial U}/{\partial x},$ whenever $U$ is sufficiently regular, we have that
$$
T_{1/2}T_{1/2}\varphi=\frac{d\varphi}{dt}.
$$
That is $T_{1/2}$ acts like an half derivative, indeed 
$$
\frac{\partial}{\partial x}\frac{\partial U}{\partial x}(x,t)=\frac{\partial U}{\partial t}(x,t) \underset{x \rightarrow 0^+} \longrightarrow  \frac{d\varphi(t)}{dt}.
$$
The solution of the problem (\ref{pro1}) under the reasonable assumptions that $\varphi$ is bounded and H{\"o}lder continuous, is explicitly known (check \cite{tichonov}, Chapter 3.3) to be
\begin{equation*}\begin{split} U(x,t)=&\, c {x} \int_{-\infty}^t \displaystyle e^{-\frac{x^2}{4(t-\tau)}}{(t-\tau)^{-\frac{3}{2}}}\varphi(\tau)\, d\tau\\
= &\,c {x}  \int_{0}^{\infty} \displaystyle e^{-\frac{x^2}{4\tau}}{\tau^{-\frac{3}{2}}}\varphi (t-\tau)\, d\tau,
\end{split}
\end{equation*}
where the last line is obtained with a change of variables. Using $\displaystyle t ={x^2}/{(4\tau)}$ and the integral definition of the Gamma function 
 we have that 
\[ \int_0^\infty x e^{-\frac{x^2}{4\tau}} \tau ^{-\frac{3}2} \, d\tau = 2 \int_0^\infty e^{-t} t^{-\frac{1}2}\, dt = 2 \Gamma\left(\frac{1}{2}\right).\]
As a consequence,
$$ \frac{U(x,t)-U(0,t)}{x}= c\int_{0}^{\infty}e^{-\frac{x^2}{4\tau}}{\tau^{-\frac{3}{2}}}\left(\varphi(t-\tau)-\varphi(t)\right)d\tau,$$
choosing $c$ that takes into account the right normalization. In addition, by passing to the limit, we obtain
\begin{equation*}
-\lim_{x\to 0^+}\frac{U(x,t)-U(0,t)}{x}=c\int_{0}^{\infty}\frac{\varphi(t)-\varphi(t-\tau)}{\tau^{\frac{3}{2}}}d\tau.
\end{equation*}
Hence, with the right choice of the constant, we get exactly $\bold{D}^{1/2}\varphi$ 
 i.e. the Marchaud derivative  of order $1/2$ of $\varphi.$

In \cite{Torreault}, and independently also in \cite{BF}, has been proved the following result.
 
\begin{teo} \label{teo:mainstat}
Let $s\in (0,1) $ and $\bar{\gamma}\in (s,1]$ be fixed. Let $\varphi \in C^{\bar{\gamma}}(\mathbb{R})$ 
be a bounded function 
and let $U\colon [0,\infty)\times \mathbb{R}\to \mathbb{R}$ be a solution of the problem 
\begin{equation}\label{prob1}
\left\{
\begin{split}
	&\frac{\partial U}{\partial t} (x,t)= \frac{1-2s}{x} \frac{\partial U }{\partial x}(x,t)+ \frac{\partial^2 U }{\partial x^2}(x,t), & &    (x,t)\in(0,\infty)\times \mathbb{R}\\
	 & U(0,t)=\varphi(t),  &&t\in\mathbb{R}\\
	  & \lim_{x \to +\infty} U(x,t)=0. & &
\end{split}
\right.
\end{equation}
Then $U$ defines the extension operator for $\varphi$, such that
\begin{equation*}
\bold{D}^s \varphi(t)=-\lim_{x\to 0^+} c_s x^{-2s}(U(x,t)- \varphi(t)),\quad \mbox{ where }\quad c_s= 4^s\Gamma(s).
\end{equation*}
 \end{teo}

An interesting application that follows from this extension procedure is a Harnack inequality for Marchaud stationary functions in an interval $J \subseteq {\mathbb{R}},$ namely for functions that satisfy  $\bold{D}^s \varphi=0$  in  $J.$ This fact is not obvious, indeed the set of functions determined by fractional-stationary functions (on an interval) is nontrivial, see e.g. \cite{Bucur}.  
\begin{teo}\label{teoHarn}
Let $s\in (0,1)$. There exists a positive constant $\gamma$ such that, if $\bold{D}^s\varphi=0$ in an interval $J\subseteq \mathbb{R}$ and $\varphi\geq 0$ in $\mathbb{R}$, then  
\begin{equation*}\label{Harnack_intro}
\sup_{[t_0-\frac{3}{4}\delta,t_0-\frac{1}{4}\delta]}\varphi\leq \gamma \inf_{[t_0+\frac{3}{4}\delta,t_0+\delta]}\varphi
\end{equation*}
for every $t_0\in \mathbb{R}$ and for every  $\delta >0$ such that $[t_0-\delta,t_0+\delta]\subset J$.
\end{teo}

The previous result can be deduced from the Harnack inequality proved in \cite{ChiSe} for some degenerate parabolic operators (see also \cite{FKS} for the elliptic setting) that however are local operators. In particular, the constant $\gamma$ used in  Theorem \ref{teoHarn} is the same that appears in the parabolic case in \cite{ChiSe}. Concerning 
the Harnack inequality for the Riemann-Liouville fractional derivative we point out also \cite{Z1} and \cite{Z2}. In concluding this section we also remark that, as  far as in the case of the fractional Laplace case, the result is true if $\phi\geq 0$ in all of $\mathbb{R},$ see \cite{Kassmann} for a counterexample for the fractional Laplace case.

We end this section remarking that
concerning the numerical computation of the fractional operators, there exist many contributions. Among them we point out \cite{DFFD} and the recent handbook \cite{LZ}.

\section{ Relationship between Marchaud derivative and the fractional Laplace operator}\label{rel_Mar_Lapl}
In the end, we discuss here some relationships between Marchaud derivative and fractional Laplace operators. An application to this approach can be find in \cite{Fer_H} in the first Heisenberg group case. By the way we would like to point out that recently a major and renewed attention to fractional calculus and operators similar to Marchaud derivative has been testified by the application described in \cite{ACV}. 

In order to explain how fractional Laplace and Marchaud derivative are linked, we fix our attention to the case $0<\alpha<1$
by considering 
$$
{\bf D}^\alpha_{\pm, \xi} f(x)=\frac{\alpha}{\Gamma(1-\alpha)}\int_0^{+\infty}\frac{f(x)-f(x\mp t\xi )}{t^{1+\alpha}}dt
$$
where $\xi\in \mathbb{S}^{n-1}$ and $f\in\mathcal{S}(\mathbb{R}^n).$
For clarity, we define a new operator as it follows: for every $f\in \mathcal{S}(\mathbb{R}^n),$
$$
\mathcal{M}_{\frac{\alpha}{2}}f(x)=\int_{\partial B_1(0)}{\bf D}^\alpha_{\xi} f(x)d\mathcal{H}^{n-1}(\xi).
$$
Then switching the order of integration we get
\begin{equation}
\begin{split}
&\int_{\partial B_1(0)}{\bf D}^\alpha_{\xi} f(x)d\mathcal{H}^{n-1}(\xi)=\frac{\alpha}{\Gamma(1-\alpha)}\int_{0}^{+\infty}\left(\int_{\partial B_1(0)}\frac{f(x)-f(x-t\xi)}{t^{1+\alpha}}d\mathcal{H}^{n-1}(\xi)\right)dt\\
&=\frac{\alpha}{\Gamma(1-\alpha)}\int_{0}^{+\infty}\left(\int_{\partial B_t(x)}\frac{f(x)-f(y)}{|x-y|^{1+\alpha}}d\mathcal{H}^{n-1}(y)\right)\frac{dt}{t^{n-1}}\\
&=\frac{\alpha}{\Gamma(1-\alpha)}\int_{0}^{+\infty}\left(\int_{\partial B_t(x)}\frac{f(x)-f(y)}{|x-y|^{n+\alpha}}d\mathcal{H}^{n-1}(y)\right)dt=\frac{\alpha}{\Gamma(1-\alpha)}\frac{\beta_n(\alpha)}{\sin(\alpha\frac{\pi}{2})}(-\Delta)^{\frac{\alpha}{2}}f(x).
\end{split}
\end{equation}
In general, as already remarked in Lemma 26.2, \cite{SKM}, and recalling the previous Section \ref{fractional_recall} for the definition of the constants $\chi(\alpha,l),$ for  every
 $\alpha>0,$ $l\in \mathbb{N},$  $l\geq 0,$ $\alpha<l,$ denoting
$$
\mathcal{D}^{\alpha)}_\xi f(x)=\frac{1}{\chi(\alpha,l)}\int_{0}^{+\infty}\frac{\Delta^l_{\xi}f(x)}{\tau^{1+\alpha}}d\tau,
$$
we obtain, for every $f\in \mathcal{S}(\mathbb{R}^n):$
\begin{equation}
\begin{split}
(-\Delta)^{\frac{\alpha}{2}}f(x)=-\frac{\Gamma(-\alpha)\sin(\alpha\frac{\pi}{2})}{\beta_n(\alpha)}\int_{\partial B_1(0)}{\mathcal{D}}^{\alpha)}_{\xi} f(x)d\mathcal{H}^{n-1}(\xi).
\end{split}
\end{equation}


In particular, if $\alpha\in]0,1[,$ and $l=1$
$$
{\bold D}^{\alpha,+}f(t)=\frac{1}{C_{\alpha,1}}\int_{0}^{+\infty}\frac{f(t)-f(t-s)}{s^{1+\alpha}}ds,
$$
where $$
C_{\alpha,1}=\frac{\Gamma (1-\alpha)}{\alpha}$$
and
$$
{\bold D}^{\alpha,-}f(t)=\frac{1}{C_{\alpha,1}}\int_{0}^{+\infty}\frac{f(t)-f(t+s)}{s^{1+\alpha}}ds,
$$
where $$
C_{\alpha,1}=\frac{\Gamma (1-\alpha)}{\alpha}.$$

Thus,
$$
({\bold D}^{\alpha,+}+{\bold D}^{\alpha,-})f(t)=\frac{1}{C_{\alpha,1}}\int_{0}^{+\infty}\frac{2f(t)-f(t+\tau)-f(t-\tau)}{\tau^{1+\alpha}}d\tau.
$$

and for every $e\in \partial B_1(0)$ and for every $f\in \mathcal{S}(\mathbb{R}^n)$ we have:
$$
({\bold D}^{\alpha,+}_e+{\bold D}^{\alpha,-}_e)f(x)=\frac{1}{C_{\alpha,1}}\int_{0}^{+\infty}\frac{2f(x)-f(x+e\tau)-f(x-e\tau)}{\tau^{1+\alpha}}d\tau.
$$
As a consequence, integrating on $\partial B_1(0),$ we obtain, as we already remarked in one variable only:
\begin{equation*}
\begin{split}
&\int_{\partial B_1(0)}({\bold D}^{\alpha,+}_e+{\bold D}^{\alpha,-}_e)f(x)d\mathcal{H}^{n-1}(e)\\
&=\frac{1}{C_{\alpha,1}}\int_{\partial B_1(0)}\left(\int_{0}^{+\infty}\frac{2f(x)-f(x+e\tau)-f(x-e\tau)}{\tau^{1+\alpha}}d\tau\right)d\mathcal{H}^{n-1}(e)\\
&=\frac{1}{C_{\alpha,1}}\int_{0}^{+\infty}\left(\int_{\partial B_1(0)}\frac{2f(x)-f(x+e\tau)-f(x-e\tau)}{\tau^{1+\alpha}}d\mathcal{H}^{n-1}(e)\right)d\tau\\
&=\frac{1}{C_{\alpha,1}}\int_{0}^{+\infty}\left(\int_{\partial B_\tau(0)}\frac{2f(x)-f(x+\xi)-f(x-\xi)}{\tau^{n+\alpha}}d\mathcal{H}^{n-1}(\xi)\right)d\tau\\
&=\frac{1}{C_{\alpha,1}}\int_{0}^{+\infty}\left(\int_{\partial B_\tau(0)}\frac{2f(x)-f(x+\xi)-f(x-\xi)}{|x-\xi|^{n+\alpha}}
d\mathcal{H}^{n-1}(\xi)\right)d\tau\\
&=\frac{1}{C_{\alpha,1}}\int_{\mathbb{R}^n}\frac{2f(x)-f(x+\xi)-f(x-\xi)}{|x-\xi|^{n+\alpha}}d\xi=\frac{2}{C_{\alpha,1}c(\frac{\alpha}{2},n)}(-\Delta)^{\frac{\alpha}{2}}f(x).
\end{split}
\end{equation*}

\end{document}